# A METHOD OF SOLVING A DIOPHANTINE EQUATION OF SECOND DEGREE WITH N VARIABLES

Florentin Smarandache

University of New Mexico

200 College Road

Gallup, NM 87301, USA

**ABSTRACT**. First, we consider the equation

(1)  $ax^2 - by^2 + c = 0$, with $a,b \in \mathbb{N}^*$ and $c \in \mathbb{Z}^*$.

It is a generalization of Pell's equation: $x^2 - Dy^2 = 1$. Here, we show that: if the equation has an integer solution and $a \cdot b$ is not a perfect square, then (1) has infinitely many integer solutions; in this case we find a closed expression for $(x_n, y_n)$, the general positive integer solution, by an original method. More, we generalize it for a Diophantine equation of second degree and with n variables of the form:

$$\sum_{i=1}^{n} a_i x_i^2 = b, \text{ with all } a_i, b \in \mathbb{Z}, n \geq 2.$$

**INTRODUCTION.**

If a·b = k$^2$ is a perfect square (k∈N) the equation (1) has at most a finite number of integer solutions, because (1) becomes:

(2)  (ax - ky) (ax + ky) = - ac.

If (a, b) does not divide c, the Diophantine equation has no solution.

**METHOD OF SOLVING.**

Suppose (1) has many integer solutions. Let $(x_0, y_0)$, $(x_1, y_1)$ be the smallest positive integer solutions for (1), with $0 \leq x_0 < x_1$. We construct the recurrent sequences:

(3) $\begin{cases} x_{n+1} = \alpha x_n + \beta y_n \\ y_{n+1} = \gamma x_n + \delta y_n \end{cases}$

setting the condition that (3) verifies (1).  It results in:

$a\alpha\beta = b\gamma\delta$ (4)

$a\alpha^2 - b\gamma^2 = a$ (5)

$a\beta^2 - b\delta^2 = -b$ (6)

having the unknowns $\alpha$, $\beta$, $\gamma$, $\delta$. We pull out $a\alpha^2$ and $a\beta^2$ from (5), respectively (6), and replace them in (4) at the square; we obtain:

(7)  $a\delta^2 - b\gamma^2 = a$ .

We subtract (7) from (5) and find

(8)  $\alpha = \pm \delta$ .

Replacing (8) in (4) we obtain

(9)  $\beta = \pm \dfrac{b}{a} - \gamma$ .

Afterwards, replacing (8) in (5), and (9) in (6), we find the same equation:

(10)  $a\alpha^2 - b\gamma^2 = a$.

Because we work with positive solutions only, we take:

$$\begin{cases} x_{n+1} = \alpha_0 x_n + (b/a)\gamma_0 y_n \\ y_{n+1} = \gamma_0 x_n + \alpha_0 y_n \end{cases},$$

where $(\alpha_0, \gamma_0)$ is the smallest positive integer solution of (10) such that $\alpha_0 \gamma_0 \neq 0$. Let the 2x2 matrix be:

$$A = \begin{pmatrix} \alpha_0 & (b/a)\gamma_0 \\ \gamma_0 & \alpha_0 \end{pmatrix} \in M_2(Z).$$

Of course, if $(x', y')$ is an integer solution for (1), then $A \cdot \begin{pmatrix} x_0 \\ y_0 \end{pmatrix}, A^{-1} \cdot \begin{pmatrix} x_0 \\ y_0 \end{pmatrix}$ is another one, where $A^{-1}$ is the inverse matrix of A, i.e., $A^{-1} \cdot A = A \cdot A^{-1} = I$ (unit matrix). Hence, if (1) has an integer solution, it has infinitely many (clearly $A^{-1} \in M_2(Z)$).

The <u>general positive integer solution</u> of the equation

(1)  is

$(x'_n, y'_n) = (|x_n|, |y_n|)$, with

$$(\text{GS}_1) \quad \begin{pmatrix} x_n \\ y_n \end{pmatrix} = A^n \cdot \begin{pmatrix} x_0 \\ y_0 \end{pmatrix}, \text{ for all } n \in \mathbb{Z},$$

where by convention $A^0 = I$ and $A^{-k} = A^{-1} \cdot \ldots \cdot A^{-1}$ of k times.
In the problems it is better to write (GS) as:

$$\begin{pmatrix} x_n' \\ y_n' \end{pmatrix} = A^n \cdot \begin{pmatrix} x_0 \\ y_0 \end{pmatrix},$$

$n \in \mathbb{N}$, and

$$(\text{GS}_2) \quad \begin{pmatrix} x_n'' \\ y_n'' \end{pmatrix} = A^n \begin{pmatrix} x_1 \\ y_1 \end{pmatrix}, \quad n \in \mathbb{N}*.$$

We prove by *reductio ad absurdum* that ($\text{GS}_2$) is a <u>general positive integer solution</u> for (1).

Let (u, v) be a positive integer particular solution for (1). If $\exists k_0 \in \mathbb{N}$: $(u, v) = A^{k_0} \cdot \begin{pmatrix} x_0 \\ y_0 \end{pmatrix}$, or

$\exists k_1 \in \mathbb{N}$: $(u, v) = A^{k_1} \cdot \begin{pmatrix} x_1 \\ y_1 \end{pmatrix}$, then $(u, v) \in (\text{GS}_2)$.

Contrarily to this, we calculate $(u_{i+1}, v_{i+1}) = A^{-1} \cdot \begin{pmatrix} u_i \\ v_i \end{pmatrix}$ for

i = 0, 1, 2, ..., where $u_0 = u$, $v_0 = v$. Clearly $u_{i+1} < u_i$

for all i. After a certain rank, $i_0$, it is found that

$x_0 < u_{i_0} < x_1$ or $0 < u_{i_0} < x_0$, but that is absurd.

It is clear we can put

$$(\text{GS}_3) \quad \begin{pmatrix} x_n \\ y_n \end{pmatrix} = A^n \cdot \begin{pmatrix} x_0 \\ \varepsilon y_0 \end{pmatrix}, \quad n \in \mathbb{N}, \text{ where } \varepsilon = \pm 1.$$

We have now to transform the general solution ($GS_3$) into a closed expression. Let $\lambda$ be a real number. $Det(A - \lambda \cdot I) = 0$ involves the solutions $\lambda_{1,2}$ and the proper vectors $v_{1,2}$ (i.e., $Av_i = \lambda_i v_i$, $i \in \{1,2\}$). Note $P = \begin{pmatrix} v_1 \\ v_2 \end{pmatrix}^t \in M_2(R)$.

Then $P^{-1}AP = \begin{pmatrix} \lambda_1 & 0 \\ 0 & \lambda_2 \end{pmatrix}$, whence $A^n = P \cdot \begin{pmatrix} (\lambda_1)^{\wedge}n & 0 \\ 0 & (\lambda_2)^{\wedge}n \end{pmatrix} \cdot P^{-1}$, and,

replacing it in ($GS_3$) and doing the calculation, we find a closed expression for ($GS_3$).

**EXAMPLES.**

1. For the Diophantine equation $2x^2 - 3y^2 = 5$ we obtain:

$$\begin{pmatrix} x_n \\ y_n \end{pmatrix} = \begin{pmatrix} 5 & 6 \\ 4 & 5 \end{pmatrix}^n \cdot \begin{pmatrix} 2 \\ \varepsilon \end{pmatrix}, \quad n \in N,$$

and $\lambda_{1,2} = 5 \pm 2\sqrt{6}$, $v_{1,2} = (\sqrt{6}, \pm 2)$, whence a closed expression for $x_n$ and $y_n$:

$$x_n = \frac{4+\varepsilon\sqrt{6}}{4}(5+2\sqrt{6})^n + \frac{4-\varepsilon\sqrt{6}}{4}(5-2\sqrt{6})^n$$

$$y_n = \frac{3\varepsilon+2\sqrt{6}}{6}(5+2\sqrt{6})^n + \frac{3\varepsilon-2\sqrt{6}}{6}(5-2\sqrt{6})^n,$$

for all n∈N.

2. For the equation $x^2 - 3y^2 - 4 = 0$ the general solution in positive integers is:

$$x_n = (2+\sqrt{3})^n + (2-\sqrt{3})^n$$

$$y_n = \frac{1}{\sqrt{3}}[(2+\sqrt{3})^n - (2-\sqrt{3})^n]$$

for all n∈N, that is (2, 0), 4, 2), (14, 8), (52, 30), ... .

**EXERCISES FOR READERS.**

Solve the Diophantine equations:

3. $x^2 - 12y^2 + 3 = 0$.

Remark:
$$\begin{pmatrix} x_n \\ y_n \end{pmatrix} = \begin{pmatrix} 7 & 24 \\ 2 & 7 \end{pmatrix}^n \cdot \begin{pmatrix} 3 \\ \varepsilon \end{pmatrix} = ?, \; n \in N.$$

4. $x^2 - 6y^2 - 10 = 0$.

Remark:
$$\begin{pmatrix} x_n \\ y_n \end{pmatrix} = \begin{pmatrix} 5 & 12 \\ 2 & 5 \end{pmatrix}^n \cdot \begin{pmatrix} 4 \\ \varepsilon \end{pmatrix} = ?, \; n \in N.$$

5. $x^2 - 12y^2 + 9 = 0$.

Remark:
$$\begin{pmatrix} x_n \\ y_n \end{pmatrix} = \begin{pmatrix} 7 & 24 \\ 2 & 7 \end{pmatrix}^n \cdot \begin{pmatrix} 3 \\ 0 \end{pmatrix} = ?, \; n \in N.$$

6.  $14x^2 - 3y^2 - 18 = 0$.

**GENERALIZATIONS.**

If $f(x, y) = 0$ is a Diophantine equation of second degree with two unknowns, by linear transformations it becomes:

(12) $ax^2 + by^2 + c = 0$, with $a, b, c \in Z$.

If $a \cdot b \geq 0$ the equation has at most a finite number of integer solutions which can be found by attempts.

It is easier to present an example:

1.  The Diophantine equation:

(13) $9x^2 + 6xy - 13y^2 - 6x - 16y + 20 = 0$

becomes:

(14) $2u^2 - 7v^2 + 45 = 0$, where

(15) $u = 3x + y - 1$ and $v = 2y + 1$.

We solve (14). Thus:

(16) $\begin{aligned} u_{n+1} &= 15u_n + 28v_n \\ v_{n+1} &= 8u_n + 15v_n \end{aligned}$, $n \in N$, with $(u_0, v_0) = (3, 3\varepsilon)$.

<u>First Solution.</u>

By induction we prove that: for all $n \in N$ we have: $v_n$ is odd, and $u_n$ as well as $v_n$ are multiples of 3. Clearly

$v_0 = 3\varepsilon$, $u_0 = 3$. For n + 1 we have: $v_{n+1} = 8u_n + 15v_n =$

= even + odd = odd, and of course $u_{n+1}$, $v_{n+1}$ are multiples of 3 because $u_n$, $v_n$ are multiples of 3, too. Hence, there exists $x_n$, $y_n$ in positive integers for all n∈N:

$$\begin{aligned} x_n &= (2u_n - v_n + 3)/6 \\ y_n &= (v_n - 1)/2 \end{aligned} \quad (17)$$

(from (15)). Now we find the ($GS_3$) for (14) as closed expression, and by means of (17) it results the general integer solution of the equation (13).

Second Solution.

Another expression of the ($GS_3$) for (13) we obtain if we transform (15) as: $u_n = 3x_n + y_n - 1$ and $v_n = 2y_n + 1$,

for all n∈N. Whence, using (16) and doing the calculation,

we find:

$$\begin{aligned} x_{n+1} &= 11x_n + \frac{52}{3} y_n + \frac{11}{3} \\ y_{n+1} &= 12x_n + 19y_n + 3, \; n \in N, \end{aligned} \quad (18)$$

with $(x_0, y_0) = (1, 1)$ or $(2, -2)$

(two infinitudes of integer solutions).

Let $A = \begin{pmatrix} 11 & 52/3 & 11/3 \\ 12 & 19 & 3 \\ 0 & 0 & 1 \end{pmatrix}$, then $\begin{pmatrix} x_n \\ y_n \\ 1 \end{pmatrix} = A^n \cdot \begin{pmatrix} 1 \\ 1 \\ 1 \end{pmatrix}$ or

$$\begin{pmatrix} x_n \\ y_n \\ 1 \end{pmatrix} = A^n \cdot \begin{pmatrix} 2 \\ -2 \\ 1 \end{pmatrix}, \text{ always } n \in N; \quad (19).$$

From (18) we have always $y_{n+1} \equiv y_n \equiv \ldots \equiv y_0 \equiv 1 \pmod{3}$, hence always $x_n \in Z$. Of course (19) and (17) are equivalent as general integer solution for (13). [The reader can calculate $A^n$ (by the same method liable to the start of this note) and find a closed expression for (19).]

**More General.**

This method can be generalized for the Diophantine equations of the form:

$$(20) \quad \sum_{i=1}^{n} a_i x_i^2 = b, \text{ with all } a_i, b \in Z, n \geq 2.$$

If $a_i \cdot a_j \geq 0$, $1 \leq i < j \leq n$, is for all pairs $(i, j)$, equation (20) has at most a finite number of integer solutions.

Now, we suppose $\exists i_0, j_0 \in \{1, \ldots, n\}$ for which $a_{i_0} \cdot a_{j_0} < 0$ (the equation presents at least a variation of sign). Analogously, for $n \in N$, we define the recurrent sequences:

$$(21) \quad x_h^{(n+1)} = \sum_{i=1}^{n} \alpha_{ih} x_i^{(n)}, \quad 1 \leq h \leq n,$$

considering $(x_1^0, \ldots, x_n^0)$ the smallest positive integer solution of (20). One replaces (21) in (20), one identifies the coefficients and one looks for the $n^2$ unknowns $\alpha_{ih}$, $1 \leq i, h \leq n$. (This calculation is very intricate, but it can be done by means of a computer.) The method goes on similarly, but the calculation becomes more and more intricate, for example to calculate $A^n$. [The reader will be able to try his/her forces for the Diophantine equation $ax^2 + by^2 - cz^2 + d = 0$, with a, b, c∈N* and d∈Z.]